\documentclass[a4paper]{amsart}
\usepackage{amssymb, amsmath, amsthm}

\theoremstyle{plain}
\newtheorem{thm}[subsection]{Theorem}

\newtheorem{lem}[subsubsection]{Lemma}

\theoremstyle{definition}
\newtheorem{deft}[subsubsection]{Definition}

\theoremstyle{remark}
\newtheorem{rem}[subsubsection]{Remark}
\newtheorem*{prob}{Problem}
\newtheorem{exam}[subsubsection]{Example}

\DeclareMathOperator{\dist}{dist}
\DeclareMathOperator{\codim}{codim}
\DeclareMathOperator{\Diag}{Diag}

\newcommand{\imax}{\underline{m}}
\newcommand{\loja}{\mathfrak{L}}

\title[Infinite determinacy for non-isolated singularities]{Infinite determinacy on a closed set for smooth germs with non-isolated singularities}

\author{Vincent Thilliez}

\address{Math\'ematiques - B\^atiment M2\\
Universit\'e des Sciences et Technologies de Lille\\
F-59655 Villeneuve d'Ascq Cedex, France}

\email{thilliez@math.univ-lille1.fr}


\begin{document}

\begin{abstract}
We give necessary and sufficient conditions of infinite determinacy for sm\-ooth func\-tion germs whose critical locus contains a
given set. This set is assumed to be the zero variety $ X $ of some analytic map-germ having maximal rank on a dense subset of $ X $. We obtain a result in terms of \L ojasiewicz estimates which extends, in particular, previous works by Sun \& Wilson on line singularities, and by Grandjean on singularities of codimension $ 1 $ or $ 2 $.    
\end{abstract}

\maketitle

\section*{Introduction}

Infinite determinacy is a way to express the stability of smooth function-germs under flat perturbations. Let $ \mathcal{E}_n $ denote the ring of $ C^\infty $ 
function-germs at the origin in $ \mathbb{R}^n $ and $ \imax $ its maximal ideal. Consider the ideal
of flat germs $ \imax^\infty=\bigcap_{k\geq 0}\imax^k $. An
element $ f $ of $ \mathcal{E}_n $ is said to be \emph{infinitely determined} if, for
any element $ u $ of $ \imax^\infty $, there exists a germ $ \Phi $ of
$ C^\infty $-diffeomorphism at the origin such that $ f+u=f\circ\Phi $. In abbreviated form, this can be written
\begin{equation}\label{deter00}
f+\imax^\infty\subset f\circ\mathcal{R}, 
\end{equation}
where $ \mathcal{R} $ denotes the group of
germs at the origin of $ C^\infty $-diffeomorphisms of $ \mathbb{R}^n $.
In what follows, for any subset $ \mathcal{R}' $ of $ \mathcal{R} $, we shall always use the notation $ f\circ\mathcal{R}' $ to denote the set
$ \{f\circ \Phi\, ;\, \Phi\in \mathcal{R}'\} $.

It is known (see \cite{Nguyen}, \cite{Wilson} or part II of \cite{Wall}) that  \eqref{deter00} holds  
if and only if
\begin{equation}\label{jac00}
\imax^\infty\subset J_f,
\end{equation}
where $ J_f $ denotes the Jacobian ideal $ \big(\frac{\partial f}{\partial x_1},\ldots \frac{\partial f}{\partial x_n}\big)\mathcal{E}_n $.  
Notice that condition \eqref{jac00} implies that $ \nabla f $ has at most an isolated singularity at $ 0 $. 

The case of of non-isolated singularities is much less understood. Typically, one considers germs with a prescribed critical locus $ X $, or germs belonging to a given ideal $ I $. If $ X $ is a line, necessary and sufficient conditions for the corresponding version of infinite determinacy have been stated first by Sun \& Wilson \cite{Sun-Wilson}. The main theorem of \cite{Sun-Wilson} appears as a nice extension of the pioneering work of Siersma \cite{Siersma} on finite determinacy for line singularities. Unfortunately, the proof of the crucial lemma 4.4 in the paper
of Sun \& Wilson is not correct, and there is currently no known correction  \cite{Wilson1}. 
Among other applications, we shall provide here a complete proof of a slightly modified version of their main statement.
To be more precise, a part of the results in \cite{Sun-Wilson} can be described as follows. The set $ X $ is the $ x_n $-axis and $ I $ is the ideal generated  in $ \mathcal{E}_n $ by the $ n-1 $ first coordinate functions
$ x_1,\ldots,x_{n-1} $. It is easy to see that the critical locus of an element $ f $ of $ \mathcal{E}_n $ contains $ X $ (as a set germ) if and only if $ f $
belongs to  $ I^2 $, and that the set of flat elements of $ I^2 $ is precisely the ideal 
$ \imax^\infty I^2 $. Accordingly, an element $ f $ of $ I^2 $ is said to be \emph{infinitely determined relatively to} $ I^2 $ if
\begin{equation}\label{deterX0}
f+\imax^\infty I^2 \subset f\circ\mathcal{R}^X, 
\end{equation}
where $ \mathcal{R}^X $ denotes the subgroup of $ \mathcal{R} $ given by the elements 
which preserve $ X $ (or, equivalently in this situation, which preserve $ I $). 
It is stated in \cite{Sun-Wilson} that \eqref{deterX0} holds if and only if
\begin{equation}\label{jacX0}
\imax^\infty I\subset J_f. 
\end{equation}
The ``only if'' part relies on another characterization of \eqref{deterX0} 
involving \L ojasiew\-icz inequalities for $ \nabla f $ and for a suitable partial
Hessian of $ f $, in the spirit of the classical formulation of \eqref{jac00} in terms of \L ojasiewicz estimates for $ \nabla f $. This is the so-called \emph{real isolated line singularity} condition in \cite{Sun-Wilson}. In order to show that
\eqref{deterX0} implies this condition, some information on the non-degeneracy of partial Hessians along the critical locus is required, and this is where a gap occurs\footnote{The computation of $ \det A'_m $ in lemma 4.4 of \cite{Sun-Wilson} is erroneous because the addition of an extra term
$ \frac{\det A_m}{\det B}y_my_l $ to the germ $ f $ modifies the Hessian symmetrically, and not only the coefficient located at the $ m $-th row and $ l $-th column.} in \cite{Sun-Wilson}. 
Our main result implies, as a particular case, that the real isolated line singularity condition characterizes a determinacy property which, instead of $ \mathcal{R}^X $, involves the (smaller) subgroup $\mathcal{R}^X_{\mathrm{fix}} $ of diffeomorphisms preserving $ X $ \emph{pointwise}. 

Although quite attractive, the real isolated line singularity condition is not easy to extend to wider classes of sets $ X $. A first step towards more general situations has been made by
Grandjean \cite{Grandjean2}, who studied the case of ideals $ I $ generated by real-analytic
germs $ \psi_1,\ldots,\psi_p $ defining a complete intersection variety $ X $ with at most an isolated singular point at the origin. It is shown in \cite{Grandjean2} that the implication
\eqref{deterX0}$\Longrightarrow$\eqref{jacX0} essentially still holds in this setting, under the additional assumption $ \codim X\leq 2 $. Notice that the proof of Grandjean is not affected by the defective lemma in \cite{Sun-Wilson}, since working in codimension $ 1 $ or $ 2 $ makes it possible to bypass the argument. On the other hand, \cite{Grandjean2} adds an \textit{a priori} assumption on the partial Hessians, namely non-degeneracy on $ X\setminus\{0\} $. As an other particular case of our main result, this extra requirement can be suppressed, provided, once again, $ \mathcal{R}^X $ is replaced by $ \mathcal{R}^X_{\mathrm{fix}} $. 

Beyond these adjustments, the main purpose of the present paper is to extend the aforementioned results by removing
the restrictive assumption on the codimension of $ X $ and, to some extent, on its singular part. In order to achieve this latter point, we involve flatness on a given closed set $ Y $, maybe larger than the single point $ 0 $. This viewpoint has been used recently in \cite{Kushner-Terra Leme} to extend the equivalence \eqref{deter00}$\Longleftrightarrow $\eqref{jac00}: indeed, 
if $ \imax_Y^\infty $
denotes the ideal of germs in $ \mathcal{E}_n $ which are flat on $ Y $,
theorem 36 of \cite{Kushner-Terra Leme} asserts that an element $ f $ of $ \mathcal{E}_n $ satisfies $ f+\imax_Y^\infty\subset f\circ\mathcal{R} $ if and only if $\imax_Y^\infty\subset J_f $. 
In order to generalize \eqref{deterX0}$\Longleftrightarrow $\eqref{jacX0} in the same way, we shall investigate here the relationship between the conditions $ f+\imax_Y^\infty I^2\subset f\circ\mathcal{R}' $ (for some suitable subgroup $ \mathcal{R}' $ of $\mathcal{R} $) and $\imax_Y^\infty I\subset J_f $.
This will lead us to the statement of theorem \ref{main}. The geometric requirements are described precisely in subsection \ref{singuls} below: roughly speaking, the regular part of $ X $ has to be dense in $ X $, and its singular part has to be contained in $ Y $. There is no restriction on $ \codim X $. Of course, the key fact is that we are able to obtain a characterization of infinite determinacy in terms of \L ojasiewicz estimates, quite in the same spirit as the real isolated line singularity condition of Sun \& Wilson.

\section{Definitions and Technical Framework}\label{basic}

\subsection{Notations} For any multi-index $ J=(j_1,\ldots,j_n) $ in $ \mathbb{N}^n $, of length $ j=j_1+\ldots+j_n $, we 
put $ D^J =\partial^j/\partial x_1^{j_1}\cdots\partial x_n^{j_n} $. 
For any finite family $ \lambda_1,\ldots,\lambda_p $ of real numbers, we denote by
$ \Diag(\lambda_1,\ldots,\lambda_p) $ the $ p\times p $ diagonal matrix whose $ j $-th diagonal coefficient is $ \lambda_j $ for $ j=1,\ldots, p$.
The differential of a map-germ $ G\, :\, (\mathbb{R}^n,0)\longrightarrow
(\mathbb{R}^m,0) $ will be denoted by $ dG $ and its jacobian matrix by $ G' $. 

\begin{deft}
Let $ Z $ be a germ of subset of $ \mathbb{R}^n $ at the origin. We 
denote
by $ \mathcal{R}^Z $ (resp. $ \mathcal{R}^Z_{\mathrm{fix}} $) the set of elements $ \Phi $ of $ \mathcal{R} $ which preserve $ Z $, that is $ \Phi(Z)\subset Z $ (resp.   which preserve $ Z $ pointwise, that is $ \Phi(x)=x $ for any $ x\in Z $).
\end{deft}

We recall, from the introduction, that $ \mathcal{R} $ denotes the
group of germs of $ C^\infty $-diffeomorph\-isms of $ (\mathbb{R}^n,0) $. 
The set $ \mathcal{R}^Z_{\mathrm{fix}} $ is obviously a subgroup of
$ \mathcal{R} $, whereas $ \mathcal{R}^Z $ may not be such.

\begin{deft} Let $ V $ and $ W $ be two germs of subsets at the origin in
$ \mathbb{R}^n $. We say that a function $ g\, :\, V\longrightarrow \mathbb{R} $ satisfies the \L ojasiewicz inequality $ \loja(V,W) $ if
there exist constants $ C>0 $ and $ \alpha\geq 0 $ such that
$ \vert g(x)\vert\geq C\dist(x,W)^\alpha\text{ for any }x\in V $.
\end{deft}

\begin{rem} Throughout the article, properties holding on a given subset $ V $ of
$ \mathbb{R}^n $ are always
understood in the sense of germs, that is on a sufficiently small representative of $ V $.
\end{rem}

\subsection{Non-isolated singularities}\label{singuls}
Let $ \psi=(\psi_1,\ldots,\psi_p) $ be a real-analytic map-germ $ (\mathbb{R}^n,0)\longrightarrow (\mathbb{R}^p,0) $ with $ p\leq n $. With $ \psi $ we associate the zero variety $ X=\psi^{-1}(\{0\}) $ and the ideal $ I=(\psi_1,\ldots,\psi_p)\mathcal{E}_n $.
We define
\begin{equation*}
\Sigma=\{x\in X\, ;\, d\psi_1(x)\wedge\cdots\wedge d\psi_p(x)=0\} 
\end{equation*}
and we make the following assumption:
\begin{equation}\label{regdense}
X\setminus\Sigma\text{\, is dense in } X.
\end{equation} 

Any sufficiently general version of the \L ojasiewicz inequality \cite{Bochnak-Risler} shows that 
\begin{equation}\label{IC}
\bigl\vert d\psi_1\wedge\cdots\wedge d\psi_p\bigr\vert\,\text{ satisfies }\loja(X,\Sigma).
\end{equation}
The analyticity of $ \psi $ also ensures that $ I $ and its powers are closed \cite{Malgrange}. Using section V.4 of \cite{Tougeron2}, we have therefore
\begin{equation}\label{I2Loj}
\imax_X^\infty\subset I^2.
\end{equation}
For any germ of closed subset $ Y $ at the origin of $ \mathbb{R}^n $, we have also, by proposition V.2.3 of \cite{Tougeron2},
\begin{equation}\label{I2clos}
\imax_Y^\infty\cap I^2=\imax_Y^\infty I^2. 
\end{equation}

Let $ J_X $ denote the set of germs $ f $ satisfying $ f(0)=0 $ and whose critical locus contains $ X $. Using local coordinates on the smooth submanifold $ X\setminus\Sigma $, it is easy to show that any element $ f $ of $ J_X $ vanishes on $ X\setminus\Sigma $, hence on $ X $ by
\eqref{regdense}. Thus, if $ I_X $ denotes the ideal of elements of $ \mathcal{E}_n $ vanishing on $ X $, it follows that $ J_X $ is the so-called \emph{primitive ideal} of $ I_X $. Obviously, $ J_X $ contains $ I^2 $. It coincides with $ I^2 $ in certain situations, as in \cite{Sun-Wilson}. The following lemma provides a simple sufficient condition for this equality.

\begin{lem}
One has $ J_X=I^2 $ whenever $ I_X=I $.
\end{lem}
\begin{proof} Assume $ I_X=I $ and let $ f $ be an element of $ J_X $. Then 
we have $ f=\sum_{j=1}^p f_j\psi_j $ with $ \sum_{j=1}^pf_j(x)d\psi_j(x)=0 $ for any $ x\in X $. Since the differentials $ d\psi_j(x) $ are independent at any point $x$ of $ X\setminus\Sigma $, we derive that each $ f_j $ vanishes on $ X\setminus\Sigma $, hence on $ X $ by \eqref{regdense}. Thus each $ f_j $ belongs to $ I $ and  $ f $ belongs to $ I^2 $.
\end{proof} 

In what follows, just as in \cite{Grandjean2}, \cite{Izumiya-Matsuoka}, \cite{Siersma} and \cite{Sun-Wilson}, the study will be limited to elements $ f $ of $ I^2 $. Beside technical considerations, the preceding lemma provides a motivation for this approach. 
Notice that the difference of two elements of $ I^2 $ is flat on
$ Y $ if and only if it belongs to $ \imax_Y^\infty I^2 $, by virtue of \eqref{I2clos}.  
Thus, we are led to consider a notion of determinacy involving $ f+\imax_Y^\infty I^2 $, as it can be seen in the statement of theorem \ref{main} below.
 
\subsection{Transversal Hessians} The Hessian matrix of an element $ f $ of $ \mathcal{E}_n $ will be denoted by $ f'' $. Assuming that $ f $ belongs to $ I^2 $, we can write $ f=\sum_{1\leq i,j\leq p}f_{ij}\psi_i\psi_j $ for some suitable elements $ f_{ij} $ of $ \mathcal{E}_n $ satisfying $ f_{ij}=f_{ji} $. Of course, the $ f_{ij} $ are not unique. However, uniqueness holds in restriction to $ X $, as shown below.

\begin{lem}\label{Hessians}\label{unique}
The matrix $ H_f=(f_{ij})_{1\leq i,j\leq p} $ satisfies 
\begin{equation}\label{hess}
{}^t \bigl(\psi'(x)\bigr)H_f(x)\psi'(x)=\frac{1}{2}f''(x)\text{ for any }x\in X.
\end{equation}
In particular, for $ x\in X $, the matrix $ H_f(x) $ is fully determined by $ \psi $ and $ f $.
\end{lem}
\begin{proof} The identity \eqref{hess} follows from a direct computation.
The uniqueness statement is an obvious consequence of \eqref{hess} since $ \psi $ has maximal rank on a dense subset of $ X $, namely on $ X\setminus\Sigma $.
\end{proof}

From now on, $ H_f $ will be considered in restriction to $ X $. At any point
of the smooth submanifold $ X\setminus\Sigma $, it can be viewed as a Hessian of $ f $ with respect to transversal directions. 

\begin{deft}
For $ x\in X $, we put $ D_f(x)=\det H_f(x) $. 
\end{deft}

As shown by lemma \ref{unique}, the function $ D_f\, :\, X\longrightarrow\mathbb{R} $ depends only on $ \psi $ and $ f $. The following elementary property will play an important role in section \ref{prove}. 

\begin{lem}\label{DfoPhi}
Let $ Z $ be a subset of $ X\setminus\Sigma $, let $ \Phi $ be an element of $ \mathcal{R}^Z $ such that both $ f $ and $ f\circ\Phi $ belong to
$ I^2 $, and let $ x $ be any point of $ Z $. Then one has $ D_{f\circ\Phi}(x)=0 $ if and only if $ D_f\big(\Phi(x)\big)=0 $.
\end{lem}
\begin{proof}
Note first that $ D_f\big(\Phi(x)\big) $ makes sense since $ \Phi(Z)\subset Z\subset X\setminus\Sigma $, in particular both $ x $ and $ \Phi(x) $ belong to $ X $.
A direct computation then
yields $ (f\circ\Phi)''(x)={}^t\bigl(\Phi'(x)\bigr)f''\big(\Phi(x)\big)\Phi'(x) $.
Using lemma \ref{Hessians}, 
we get $ {}^t\bigl(\psi'(x)\bigr)H_{f\circ\Phi}(x)\psi'(x)={}^t\bigl( (\psi\circ\Phi)'(x)\bigr)H_f\big(\Phi(x)\big)(\psi\circ\Phi)'(x) $. The
result follows, since both $ \psi'(x) $ and $ (\psi\circ\Phi)'(x) $ have maximal rank.
\end{proof}

\subsection{Jacobian and Fitting ideals}\label{jacandfit} 
In the setting of \ref{singuls}, consider the map 
$ \sigma\, :\, \mathcal{E}_n^p \longrightarrow \mathcal{E}_n $ given by 
$ \sigma(f_1,\ldots,f_p)=\sum_{i=1}^pf_i\psi_i $. Since $ \mathcal{E}_n $ is flat over the ring $ \mathcal{O}_n $ of analytic function-germs, the module of smooth relations between $ \psi_1,\ldots,\psi_p $, that is $ \ker\sigma $, admits a finite system of generators $ k^1,\ldots,k^q $ belonging
to $ \mathcal{O}_n^p $. We can assume that this system includes all the trivial relations 
$ \psi_se_r-\psi_re_s $, where $ (e_1,\ldots,e_p) $ is the canonical basis of $ \mathcal{E}_n^p $ and $ r $, $s $ are integers with $ 1\leq r<s\leq p $. 

Now let $ f $ be an element of $ I^2 $ and put $ M_f=\sigma^{-1}(J_f) $, where $ J_f $ denotes the Jacobian ideal of $ f$, as in the introduction. For any 
$ j=1,\ldots,n $, it is easy to check that the element $ h^j=(h^j_1,\ldots,h^j_p) $ of $ \mathcal{E}_n^p $ defined by 
\begin{equation}\label{Lcols}
h_i^j=\sum_{k=1}^p\Bigl(2f_{ik}\frac{\partial\psi_k}{\partial x_j}+\frac{\partial f_{ik}}{\partial x_j}\psi_k\Bigr)\quad\text{for }i=1,\ldots,p
\end{equation} 
satisfies $ \sigma(h^j)=\partial f/\partial x_j $. Then, obviously, $ M_f $ is the submodule
of $ \mathcal{E}_n^p $ generated by $ h^1,\ldots,h^n,k^1,\dots,k^q $, so that we can write $ M_f=\lambda(\mathcal{E}_n^{n+q}) $, where 
$ \lambda\, :\, \mathcal{E}_n^{n+q}
\longrightarrow \mathcal{E}_n^p $ is the morphism of free modules defined by
\begin{equation*}
\lambda(\xi_1,\dots,\xi_{n+q})= \sum_{j=1}^n\xi_jh^j+\sum_{j=1}^q\xi_{n+j}k^j.
\end{equation*}
Denote by $ \Lambda $ the matrix of $ \lambda $ in the canonical bases of $ \mathcal{E}_n^{n+q}
$ and $ \mathcal{E}_n^p $ and consider, as in \cite{Grandjean2} or \cite{Sun-Wilson}, 
the ideal $ K_f $ generated in $ \mathcal{E}_n $ by the minors of order $ p $ of $ \Lambda $. 
Following section 20.2 of \cite{Eisenbud}, we see that $ K_f $ is precisely 
the Fitting ideal $ \mathsf{Fitt}_0\big(\mathcal{E}_n^p/M_f\big) $ as it can be seen from the free presentation
\begin{equation*}
\mathcal{E}_n^{n+q}\stackrel{\lambda}{\longrightarrow}\mathcal{E}_n^p\longrightarrow\mathcal{E}_n^p/M_f\longrightarrow 0.
\end{equation*}
In particular, $ K_f $ depends only on $ \psi $ and $ f $ and it annihilates $  \mathcal{E}_n^p/M_f $, hence the inclusion
$ K_f\mathcal{E}_n^p\subset M_f $.
Taking the image by $ \sigma $, we get therefore
\begin{equation}\label{dol}
 K_fI\subset J_f.
\end{equation}
We have now all the tools required for the statement of our result.
 
\section{A Theorem of Infinite Determinacy}\label{state}

\begin{thm}\label{main}
Let $ \psi $, $ X $ and $ I $ be defined as in subsection \ref{singuls}, and let $ Y $ be a germ of closed subset of $ \mathbb{R}^n $ at the origin satisfying
\begin{equation}\label{onY}
\Sigma\subset Y.
\end{equation}
Then for any element $ f $ of $ I^2 $, the following conditions are equivalent:
\begin{equation}\label{deter}
f+\imax_Y^\infty I^2\subset f\circ\mathcal{R}^{X\setminus Y}_{\mathrm{fix}},
\end{equation}
\begin{equation}\label{ijac}
\imax_Y^\infty I\subset J_f,
\end{equation}
\begin{equation}\label{ifit}
\imax_Y^\infty \subset K_f,
\end{equation}
\begin{equation}\label{estims}
\vert\nabla f\vert\text{ satisfies }\loja(\mathbb{R}^n,X\cup Y)\text{ and } D_f\text{ satisfies }\loja(X,Y).
\end{equation}
Moreover, if one assumes additionally 
\begin{equation}\label{onYstr}
D_f^{-1}(\{0\})\subset Y,
\end{equation}
then \eqref{deter} is equivalent to the (a priori weaker) condition
\begin{equation}\label{deterwk}
f+\imax_Y^\infty I^2\subset f\circ\mathcal{R}^{X\setminus Y}.
\end{equation}
\end{thm}

All of section \ref{prove} will be devoted to proving the theorem. We shall first complete section \ref{state} with some comments and examples.

\subsection{Comments}

\begin{rem}\label{Yis0}
The particular case $ Y=\{0\} $ corresponds to the natural notion of determinacy by the Taylor jet at a single point. In view of \eqref{onY}, our result covers this case provided $ X $ has at most an isolated singular point. Notice that we then have obviously $ \mathcal{R}^{X\setminus Y}=\mathcal{R}^X $.
\end{rem}

\begin{rem}
In this special situation, it would not be difficult to add to the statement of theorem \ref{main} another equivalent condition involving suitable substitutes for the tangent spaces to group orbits which are customary in such problems. We shall not describe this condition nor use it here (see, however, example \ref{repair} below).
\end{rem}

\subsection{Examples}

\begin{exam}[isolated critical points] From theorem \ref{main}, we can easily recover the equivalence \eqref{deter00}$\Longleftrightarrow$\eqref{jac00} of the introduction, and more generally theorem 36 of \cite{Kushner-Terra Leme}.
Indeed, take $ p=n $ and $ \psi_j(x)=x_j $ so that $ I=\imax $, $ X=\{0\} $ and $ \Sigma=\emptyset $. Then \eqref{onY} is trivial, and it is easy to check that $ \imax^\infty_Y I^2=\imax^\infty_YI=\imax^\infty_Y $. The result follows.
\end{exam}

\begin{exam}[isolated critical lines]\label{repair}
The situation studied in \cite{Sun-Wilson} corresponds to $ p=n-1 $, $ \psi_j(x)=x_j $, so that $ X $ is the $ x_n $-axis. For $ Y=\{0\} $, condition \eqref{estims} coincides with the real isolated line singularity condition of Sun \& Wilson. Taking remark \ref{Yis0} into account, we recover, as announced in the introduction, the main result of \cite{Sun-Wilson} with 
$ \mathcal{R}^X $ replaced by $ \mathcal{R}^X_{\mathrm{fix}} $. 
Explicitely, an element $ \Phi $ of $ \mathcal{R} $ belongs to $ \mathcal{R}^X $ if and only if the components $ \Phi_1,\ldots,\Phi_{n-1} $ belong to $ I $, and an element $ \Phi $ of $ \mathcal{R}^X $ belongs to $ \mathcal{R}^X_{\mathrm{fix}} $ if and only if  the component $ \Phi_n $ satisfies
$ \Phi_n(x)=x_n+\theta(x) $ with
$ \theta\in I $. The corresponding tangent spaces are respectively given by
\begin{equation*}
\mathcal{T}^Xf=\Big(\frac{\partial f}{\partial x_1},\ldots,\frac{\partial f}{\partial x_{n-1}}\Big)I+\frac{\partial f}{\partial x_n}\imax\quad\text{ and }\quad 
\mathcal{T}^X_{\mathrm{fix}}f=IJ_f.
\end{equation*}
\end{exam}

\begin{exam}[complete intersections with isolated singularities]
In the case studied in \cite{Grandjean2}, one has $ \Sigma=\{0\} $ and $ Y=\{0\} $. In this situation, \cite{Grandjean2} provides 
the implication \eqref{deterwk}$\Longrightarrow$\eqref{estims} under the additional assumptions \eqref{onYstr} and $ \codim X\leq 2 $. Thus, our theorem \ref{main} removes the restriction on the codimension $ p $ of $ X $. It also removes the assumption on 
$ D_f^{-1}(\{0\}) $ as soon as $ \mathcal{R}^X $ is replaced by $ \mathcal{R}^X_{\mathrm{fix}} $.
\end{exam}

\begin{prob} Is it possible to show that \eqref{deterwk} implies \eqref{estims} without the extra assumption \eqref{onYstr}, at least when $ Y=\{0\} $~? 
\end{prob}

\section{Proofs}\label{prove}

\subsection{Proof of \eqref{deter}$\Longrightarrow$\eqref{estims}}\label{subsestims} 
We separate the estimates for $ \nabla f $ and $ D_f $.

\subsubsection{\L ojasiewicz estimate for $ \nabla f $} We follow, with suitable modifications, the proof of the corresponding inequality in \cite{Sun-Wilson}. Assume that the
estimate does not hold. One can then find a sequence $ (x_\nu)_{\nu\geq 1} $ of points of $ \mathbb{R}^n $ converging to $ 0 $ and such that 
\begin{equation}\label{smallgrad}
\vert\nabla f(x_\nu)\vert<\dist(x_\nu,X\cup Y)^\nu\text{ for any }\nu\geq 1.
\end{equation}
Extracting a subsequence if necessary, one can assume that
$ \vert x_{\nu+1}\vert<\frac{1}{3}\dist(x_\nu,X\cup Y) $ (which implies, in particular, that
$ \dist(x_\nu,X\cup Y) $ decreases). Denote by $ B_\nu $ the open euclidean ball centered at $ x_\nu $, with radius $ \frac{1}{3}\dist(x_\nu,X\cup Y) $. Then the $ B_\nu $ are pairwise disjoint, the open set $ \mathcal{U}=\bigcup_{\nu\geq 1}B_\nu $ is contained in $ \mathbb{R}^n\setminus (X\cup Y) $ and we have $ \overline{\mathcal{U}}\cap(X\cup Y)=\{0\} $. Using a classical construction of cutoff functions (see e.g. \cite{Malgrange}, chapter I, lemma 4.2), we also have a sequence $ (\chi_\nu)_{\nu\geq 1} $ of $ C^\infty $ functions such that each $ \chi_\nu $ is supported in $ B_\nu $, identically equal to $ 1 $ in a neighborhood of the center $ x_\nu $, and satisfies, for any integer $ j\geq 0 $, any multi-index $ J $ of length $ j $ and any
point $ x $ in $ \mathbb{R}^n $, the estimate
\begin{equation}\label{trunc}
\vert D^J\chi_\nu(x) \vert\leq C(j)\dist(x_\nu,X\cup Y)^{-j}
\end{equation}
where $ C(j) $ is a constant depending only on $ j $.
Beside this, Sard's lemma ensures the existence of a regular value
$ c_\nu $ of $ f $ such that $ \vert f(x_\nu)-c_\nu\vert\leq \dist(x_\nu,X\cup Y)^\nu $. Consider the affine function $ u_\nu $ given by
$ u_\nu(x)=f(x_\nu)-c_\nu +\langle\nabla f(x_\nu),x-x_\nu\rangle $. Taking into account 
\eqref{smallgrad} and \eqref{trunc}, as well as the support condition for $ \chi_\nu $, one checks that the series
$ u=\sum_{\nu\geq 1}\chi_\nu u_\nu $ 
defines an element of $ \imax_{X\cup Y}^\infty $ such that
$ u(x_\nu)=f(x_\nu)-c_\nu $ and $ \nabla u(x_\nu)=\nabla f(x_\nu) $ for any
$ \nu\geq 1 $. Since $ \imax_{X\cup Y}^\infty=\imax_X^\infty\cap\imax_Y^\infty $, we
see, using  \eqref{I2Loj} and \eqref{I2clos}, that $ u $ belongs to
$ \imax_Y^\infty I^2 $. By \eqref{deter}, one can therefore find a germ
of diffeomorphism $ \Phi $ such that 
$ f\circ\Phi=f-u $. But this is impossible since each $ c_\nu $ is a regular value of $ f $ and a singular value of $ f-u $.

\subsubsection{\L ojasiewicz estimate for $ D_f $}\label{lojaD}
Assume that $ D_f $ does not satisfy $ \loja(X,Y) $. One can then 
find a sequence $ (y_\nu)_{\nu\geq 1} $ of points of $ X $ converging to $ 0 $ and such that
\begin{equation}\label{smallDf}
\vert D_f(y_\nu)\vert<\dist(y_\nu,Y)^\nu\quad\text{for any }\nu\geq 1. 
\end{equation}
Denote by $ \lambda_1^{(\nu)},\ldots,\lambda_p^{(\nu)} $ the eigenvalues of the symmetric matrix $ H_f(y_\nu) $, counted with multiplicities. Since we have
$ D_f(y_\nu)=\lambda_1^{(\nu)}\cdots\lambda_p^{(\nu)} $, the estimate \eqref{smallDf} implies, for each $ \nu\geq 1 $, the existence of at least one index $ i_\nu $ such that 
\begin{equation}\label{smalleig}
\big\vert \lambda_{i_\nu}^{(\nu)}\big\vert <\dist(y_\nu,Y)^{\nu/p}.
\end{equation}
Beside this, one can also find an orthogonal matrix
$ P_\nu $ such that 
\begin{equation}\label{eigenH}
P_\nu^{-1}H_f(y_\nu)P_\nu=\Diag\big(\lambda_1^{(\nu)},\ldots,\lambda_p^{(\nu)}\big).
\end{equation}
We define a $ p\times p $ symmetric matrix $ V_\nu $ by putting 
\begin{equation}\label{eigenU}
P_\nu^{-1} V_\nu P_\nu=\Diag\big(0,\ldots,0,\lambda_{i_\nu}^{(\nu)},0,\ldots,0\big),
\end{equation}
where $ \lambda_{i_\nu}^{(\nu)} $ is at the $ i_\nu $-th position. Using the fact that $ \Vert P_\nu^{-1}\Vert=\Vert P_\nu\Vert=1 $ for the euclidean norm, it is easy to see that the coefficients $ v_{ij}^{(\nu)} $ of 
$ V_\nu $ satisfy $ \big\vert v_{ij}^{(\nu)}\big\vert\leq\big\vert \lambda_{i_\nu}^{(\nu)}\big\vert $, hence, by virtue of \eqref{smalleig}, 
\begin{equation}\label{controlu}
\big\vert v_{ij}^{(\nu)}\big\vert\leq \dist(y_\nu,Y)^{\nu/p}.
\end{equation}
Now choose real numbers $ \varepsilon_1^{(\nu)},\ldots,\varepsilon_p^{(\nu)} $ satisfying, for any $ i=1,\ldots,p $,
\begin{equation}\label{epsi}
\lambda_i^{(\nu)} \neq \varepsilon_i^{(\nu)}\quad
\text{and}\quad\big\vert\varepsilon_i^{(\nu)}\big\vert<\dist(y_\nu,Y)^\nu, 
\end{equation}
then define a $ p\times p $ symmetric matrix $ W_\nu $ by putting
\begin{equation}\label{eigenV}
P_\nu^{-1}W_\nu P_\nu=\Diag\big(\varepsilon_1^{(\nu)},\ldots,\varepsilon_p^{(\nu)}\big).
\end{equation}
In view of the second condition in \eqref{epsi}, the same argument as for $ V_\nu $ shows that the coefficients $ w_{ij}^{(\nu)} $ of 
$ W_\nu $ satisfy
\begin{equation}\label{controlv}
\big\vert w_{ij}^{(\nu)}\big\vert\leq \dist(y_\nu,Y)^\nu.
\end{equation}
We can also assume that $ \vert y_{\nu+1}\vert<\frac{1}{3}\dist(y_\nu,Y) $ and consider the open ball $ C_\nu $ centered at $ y_\nu $ with radius $ \frac{1}{3}\dist(y_\nu,Y) $. With each $ C_\nu $ we associate a cutoff function $ \eta_\nu $ supported in $ C_\nu $, identically equal to $ 1 $ in a neighborhood of $ y_\nu $, and such that
\begin{equation}\label{retrunc}
\big\vert D^J\eta_\nu(x)\big\vert\leq C(j)\dist(y_\nu,Y)^{-j}
\end{equation}
for any integer $ j\geq 0 $, any multi-index $ J $ of length $ j $ and any $ x\in\mathbb{R}^n $. The open set
$ \mathcal{V}=\bigcup_{\nu\geq 1}C_\nu $ is contained in $ \mathbb{R}^n\setminus Y $ and satisfies $ \overline{\mathcal{V}}\cap Y=\{0\} $. 
By \eqref{controlu}, \eqref{controlv} and the support properties of $ \eta_\nu $, both  series
$ V=\sum_{\nu\geq 1}\eta_\nu V_\nu $ and $ W=\sum_{\nu\geq 1}\eta_\nu W_\nu $ define  symmetric $ p\times p $ matrices whose respective coefficients $ v_{ij} $ and $ w_{ij} $ all belong to $ \imax_Y^\infty $. 
Therefore, the functions $ v=\sum_{1\leq i,j\leq p}v_{ij}\psi_i\psi_j $ and
$ w=\sum_{1\leq i,j\leq p}w_{ij}\psi_i\psi_j $ belong
to $ \imax_Y^\infty I^2 $ and, by assumption, one can then find $ \Phi_v $ and $ \Phi_w $ in $ \mathcal{R}^{X\setminus Y}_{\mathrm{fix}} $ such that
\begin{equation}\label{composites} 
f-v=f\circ\Phi_v\quad \text{and}\quad f-w=f\circ\Phi_w.
\end{equation} 
Beside this, we have $ H_{f-v}(y_\nu)=H_f(y_\nu)-V(y_\nu)=H_f(y_\nu)-V_\nu $ by contruction. Therefore, \eqref{eigenH} and 
\eqref{eigenU} imply that $ 0 $ is an eigenvalue of $ H_{f-v}(y_\nu) $, hence $ D_{f-v}(y_\nu)=0 $. Thanks to \eqref{onY}, we can now use lemma \ref{DfoPhi} with
$ Z=X\setminus Y $ and $ x=y_\nu $. Thus, \eqref{composites} implies
\begin{equation}\label{Phiu}
D_f\big(\Phi_v(y_\nu)\big)=0\quad\text{for any }\nu\geq 1.
\end{equation}
We have similarly $ H_{f-w}(y_\nu)=H_f(y_\nu)-W(y_\nu)=H_f(y_\nu)-W_\nu $, hence
$ D_{f-w}(y_\nu)=\prod_{i=1}^p\big(\lambda_i^{(\nu)}-\varepsilon_i^{(\nu)}\big)\neq 0 $ by \eqref{eigenH}, \eqref{eigenV} and the first condition in \eqref{epsi}. By \eqref{composites} and lemma \ref{DfoPhi}, we get now
\begin{equation}\label{Phiv}
D_f\big(\Phi_w(y_\nu)\big)\neq0\quad\text{for any }\nu\geq 1.
\end{equation}
Recall finally that the diffeomorphisms $ \Phi_v $ and $ \Phi_w $ provided by \eqref{deter} preserve $ X\setminus Y $ pointwise, hence $ \Phi_v(y_\nu)=\Phi_w(y_\nu)=y_\nu $ and \eqref{Phiu} contradicts \eqref{Phiv}. The proof is complete.

\subsection{Proof of \eqref{estims}$\Longrightarrow$\eqref{ifit}} We follow the pattern of the corresponding proof in \cite{Sun-Wilson}, with some necessary adaptations. First, one checks
that $ K_f $ contains all the germs $ \psi_i^{p-2}\partial f/\partial x_j $ with $ 1\leq i\leq p$ and $ 1\leq j\leq n $: it suffices to compute the determinant whose first column is given by $ h^j $ and the other columns by the trivial relations $ \psi_ie_k-\psi_ke_i $ with $ k\neq i $. Therefore, $ K_f $ contains the germ
$  v=\vert\psi\vert^{2(p-2)}\vert\nabla f\vert^2 $. For $ 0<\varepsilon<1 $ and $ s\geq 1 $, put now
\begin{equation*}
V_{\varepsilon,s}=\{x\in\mathbb{R}^n\, ;\, \dist(x,X)\geq \varepsilon\dist(x,Y)^s\}\,\text{ and }\, W_{\varepsilon,s}=\mathbb{R}^n\setminus V_{\varepsilon,s}.
\end{equation*}
Since $ \vert\psi\vert^2 $ satisfies $ \loja(\mathbb{R}^n,X) $ by the classical \L ojasiewicz inequality for analytic functions, and since $ \nabla f $ satisfies $ \loja(\mathbb{R}^n,X\cup Y) $ by assumption, we see that for any given choice of $ \varepsilon $ and $ s $, the germ $ v $ satisfies $ \loja(V_{\varepsilon,s},Y) $.

Now, let 
$ \mathcal{M} $ be the set of $ p\times p $ minors of the jacobian matrix $ \psi' $.
Using \eqref{Lcols} and elementary linear algebra, one obtains, for any
$ \mu $ in $ \mathcal{M} $,
\begin{equation}\label{mudet}
2^p\mu \det(f_{ij})=a_\mu+b_\mu
\end{equation}
where $ a_\mu $ is a suitable $ p\times p $ minor of the matrix $ \Lambda $ defined in \ref{jacandfit}, hence an element of $ K_f $, and
$ b_\mu $ is
an element of $ I $.  Put
$ w=\sum_{\mu\in\mathcal{M}}a_\mu^2 $ and,
for any $ x \in (\mathbb{R}^n,0) $, denote by $ \hat{x} $ a point of $ X $ such that
$ \vert x-\hat{x}\vert=\dist(x,X) $. We have obviously   
\begin{equation}\label{eq1}
\vert w(x)-w(\hat{x})\vert\leq C\dist(x,X) 
\end{equation}
for some constant $ C>0 $. By \eqref{mudet}, we have also
\begin{equation}\label{eq2} 
w(\hat{x})=2^{2p}\bigl(D_f(\hat{x})\bigr)^2\sum_{\mu\in\mathcal{M}} \big(\mu(\hat{x})\big)^2
\end{equation}
since each $ b_\mu $ vanishes on $ X $. 
Now, by \eqref{IC} and \eqref{onY}, there exist constants $ C'>0 $ and $ \alpha\geq 1 $ such that
$ \sum_{\mu\in\mathcal{M}}\big(\mu(\hat{x})\big)^2\geq C'\dist(\hat{x},\Sigma)^\alpha 
\geq C'\dist(\hat{x},Y)^\alpha $.
Beside this, $ D_f $ satisfies $ \loja(X,Y) $ by assumption. From \eqref{eq2} and these remarks, we get
\begin{equation}\label{eq3}
\vert w(\hat{x})\vert\geq C'' \dist(\hat{x},Y)^\beta
\end{equation}
for some suitable constants $ C''>0 $ and $ \beta\geq 1 $. Assume now that $ x $
belongs to some set $ W_{\varepsilon,s} $. In this situation, we have
$ \big\vert\dist(x,Y)-\dist(\hat{x},Y)\big\vert\leq \vert x-\hat{x}\vert\leq \varepsilon \dist(x,Y)^s $, therefore $ \dist(\hat{x},Y)^\beta $ can be replaced by $ \dist(x,Y)^\beta $ in \eqref{eq3}, up to a modification of $ C'' $. From \eqref{eq1} and \eqref{eq3}, we derive that $ w $ satisfies $ \loja(W_{\epsilon,s},Y) $ for $ s>\beta $. Finally $ v+w $ is an element of $ K_f $ which satisfies
$ \loja(\mathbb{R}^n,Y) $, and proposition V.4.3 of \cite{Tougeron2} yields the conclusion. 

\subsection{Proof of \eqref{ifit}$\Longrightarrow$\eqref{ijac}} 
It is an immediate consequence of \eqref{dol}.

\subsection{Proof of \eqref{ijac}$\Longrightarrow$\eqref{deter}}\label{subsdeter}
This can be proved by a Mather-type homotopy argument, using Nakayama's lemma and tangent spaces as in \cite{Izumiya-Matsuoka} for instance. Alternatively, we use here the somewhat shorter approach based on Tougeron's implicit function theorem, see \cite{Tougeron2}, theorem III.3.2 and remark VIII.3.7.1.
Let $ u $ be an element of $ \imax_Y^\infty I^2 $ and put $ g=f+u $. For $ x $ and $ y $ in $ (\mathbb{R}^n,0) $, define $ F(x,y)=f(x+y)-g(x) $. Put $ F_i(x)=\partial F/\partial y_i(x,0) $ for $ i=1,\ldots,n $. Since
$ F_i= \partial f/\partial x_i $, it turns out that $ (F_1,\ldots,F_n)\mathcal{E}_n=J_f $. Now remark that $ \imax_Y^\infty=\imax_Y^\infty\imax_Y^\infty $ by proposition V.2.3 of \cite{Tougeron2}, hence $ \imax_Y^\infty I^2=\imax_Y^\infty(\imax_Y^\infty I)^2 $. Therefore, \eqref{ijac} implies that $ F(\,\cdot\,,0) $ belongs to $ \imax_Y^\infty J_f^2 $. Tougeron's implicit function theorem yields a map $ \varphi\,:\, (\mathbb{R}^n,0)
\longrightarrow(\mathbb{R}^n,0) $ with components in $ \imax_Y^\infty J_f $, hence in $ \imax_Y^\infty I $, such that $ F(x,\varphi(x))=0 $. Put $ \Phi(x)=x+\varphi(x) $. Clearly $ \Phi $ is a germ
of diffeomorphism at the origin, and it coincides with the identity on
$ X\cup Y $. In particular it belongs to $ \mathcal{R}^{X\setminus Y}_{\mathrm{fix}} $. This ends the proof, since we have also
$ f\circ\Phi=g $ by construction. 

\subsection{Proof of \eqref{deter}$\Longleftrightarrow$\eqref{deterwk} under the extra assumption \eqref{onYstr}}\label{weakdeter} The implication \eqref{deter}$\Longrightarrow$\eqref{deterwk} is obvious (and does not require \eqref{onYstr}). In order to prove the converse, it suffices to show that
\eqref{onYstr} together with \eqref{deterwk} imply \eqref{estims}. Now, a simple inspection of the proof of the implication \eqref{deter}$\Longrightarrow$\eqref{estims} reveals that the pointwise preservation of $ X\setminus Y $ is used only to obtain a contradiction between \eqref{Phiu} and \eqref{Phiv}. If
\eqref{deter} is replaced by \eqref{deterwk}, we obtain in the same way \eqref{Phiu} for some $ \Phi_v $ preserving $ X\setminus Y $ globally. Taking \eqref{onYstr} into account, we obtain another contradiction, since $ D_f $ has no zero on $ X\setminus Y $, hence no zero at $ \Phi_v(x_\nu) $. The result follows.

\end{document}